\documentclass[12pt]{article}
\usepackage{amsmath}
\usepackage{amssymb}
\usepackage{url}

\usepackage[T1]{fontenc}
\usepackage{mathptmx}
\usepackage{diagrams}
\diagramstyle{midshaft,PostScript=dvips,nohug}
\newarrow{Dashto}{}{.}{}{.}{>}

\newtheorem{lemma}{Lemma.}[section]
\newtheorem{prop}[lemma]{Proposition}
\newtheorem{thm}[lemma]{Theorem}
\newtheorem{cor}[lemma]{Corollary}
\newtheorem{defn}[lemma]{Definition}

\newcommand{\X}{\mathcal{X}}
\newcommand{\Y}{\mathcal{Y}}

\newcommand{\B }{\mathcal{B}}
\newcommand{\C}{\mathcal{C}}
\newcommand{\F}{\mathcal{F}}

\newcommand{\coc}{\triangleright}
\newcommand{\dcoc}{\triangleleft}

\title{Fibrations as Eilenberg-Moore algebras}
\author{Anders Kock\\
\small Dept.\ of Mathematics, University of Aarhus \normalsize}
\date{}

\begin{document}
\maketitle
\section*{Introduction}
We give an elementary account of some fundamental facts about fibered 
(or rather opfibered) categories, in terms of monads and 
2-categories. The account avoids any mention of category-valued 
functors and pseudofunctors.

I make no claim to originality;  a large part of it may  in substance be found in the  
early French literature on the subject (early 1960s -  Grothendieck, Chevalley, 
Giraud, Benabou), as reported in [Gray 1966].  The remaining  part makes use of some tools that were 
not fully available then:  monads, and their algebras, and KZ 
monads; it is extracted from   [Kock 1973, 1995] and [Street 
1974]. Some 
more historical comments are found in the  last Section.

\section{Cocartesian arrows}\label{1xx}
Given a functor $\pi :\X \to \B$. For $a:A\to B$ in $\B$, we let 
$\hom_{a}^{(\X ,\pi)}$  denote the set of arrows $x$ in $\X$ with 
$\pi (x)=a$; 
if $\pi : \X \to \B$ is clear 
from the context, we write just $\hom_{a}$. We fix a $\pi :\X \to \B$ in this Section.
 
 If 
$a$ is the identity arrow of $A$, $\hom_{a}^{\X}$ is also written 
$\hom_{A}^{\X}$; it consists of the {\em vertical} arrows of $\X$ 
over $A$ 
(relative to $\pi:\X \to \B$)\footnote{Sometimes, one needs to say
 ``$\pi$-vertical'' rather than just ``vertical'', 
namely in contexts where one  also wants to talk about ``vertical arrows'' 
meaning arrows displayed 
vertically in the graphics of a certain diagram. Often in diagrams, 
one likes to display $\pi$-vertical arrows by graphically vertical 
arrows.}. We denote by $\X_{A}$ the category 
whose objects are the objects $X\in \X$ with $\pi (X)=A$, and whose arrows are 
the vertical arrows over $A$. It is a (non-full) subcategory of 
$\X$, often called the fibre over $A$.

Given  arrows $x:X\to Y$ in  and $y:Y\to Z$ in $\X$. Let $a$ and $b$ 
denote $\pi (x)$ and $\pi (y)$, respectively, thus $x\in 
\hom_{a}(X,Y)$ and $y\in \hom_{b}(Y,Z)$. Then 
 $a$ and $b$ are composable in $\B$, and 
$x.y \in \hom_{a.b}(X,Z)$\footnote{We compose arrows in an abstract category from left to right, $a.b$ 
means ``first $a$, then $b$; whereas we compose functors between 
given categories from right to left, thus $G\circ F$ means ``first $F$, 
then  $G$''.}  
. Thus, 
for fixed $x\in \hom_{a}(X,Y)$, we have a map ``precomposition with $x$'',
\begin{equation}\label{cruc}x^{*}: \hom _{b}(Y,Z) \to \hom 
_{a.b}(X,Z).\end{equation}
The following  notion is the crucial one for the presentation here. 
\begin{defn}The arrow $x\in \hom_{a}(X,Y)$ is {\em cocartesian}\footnote{or 
op-cartesian, or, cf.\ [Johnstone 2002], supine. Note that in the classical 
definition ([Giraud 1971] p.\ 18) ``(co-)cartesian morphism'' means something 
weaker, namely as above, but with $b$ an identity arrow; an op-fibered category is 
then defined as one where there are enough of these ``weakly'' 
(co-)cartesian arrows, {\em and} where such arrows compose. In this 
case, the weak and strong notions coincide. Thus, in the 
set up of loc.cit., the notion of op-fibered category is needed prior 
to the definition of coartesian arrow, in the strong sense as given 
above. See [Borceux 1994] 8.1 for a comparison of the weak (``pre-'') and 
the strong notion.   
} if for all arrows $b$ in $\B$ 
with domain $\pi 
(Y)$ 
and all $Z \in \X_{C}$ (where $C$ denotes the codomain of $b$), 
the map (\ref{cruc}) is bijective. 
\end{defn}
Note that the only ``data'' in the definition (besides $\pi :\X \to 
\B$) is the arrow $x$. To keep track of the ``book-keeping'' involved, we 
display a diagram, in which the symbol ``$:$'' is meant to indicate 
``goes by $\pi$ to \ldots''.  The  ``data'' $a$, $A$, and $B$ are 
derived from $x$ (with $a=\pi (x)$, and $A$ and $B$ the domain and 
codomain of $a$); and $b, C$ and $Z$ are arbitrary.
\begin{equation}\label{display1}\begin{diagram} X&\rTo^{x}&Y&&Z\\
:&:&:&&:\\
A&\rTo_{a}&B&\rTo_{b}&C.
\end{diagram}\end{equation}
Another way of describing when an arrow $x$ (over $a$, say) is 
cocartesian is to say
that it has a certain (co-)universal property: for any arrow in $\X$ with 
same domain as $x$ and living over a composite arrow of the form  
$a.b$,
factorizes  as $x.y$ for  a unique $y$ over $b$. This is 
reminiscent of the (co-)universal property of a coequalizer $x$: ``any 
arrow with same domain as $x$ (and with a certain property) factors 
uniquely over $x$.''

We have the following, in case the composite $x.y$ in $\X$ is defined:
\begin{prop}\label{transx}Suppose $x$ cocartesian.  Then $y$ 
is cocartesian iff $x.y$ is cocartesian.
\end{prop}
{\bf Proof.} Straightforward verification; or see [Borceux 1994] Section 
8.1.

\medskip

We may read the 
bijections in (\ref{cruc}) as a universal property of cocartesian 
arrows $x$. Using this viewpoint, one  gets
\begin{prop}\label{vertiso} If $x:X\to Y$ and $x':X\to Y'$ are cocartesian arrows 
with $\pi (x)=\pi (x')$, then there exists a unique vertical 
isomorphism $t:Y\to Y'$ with $x.t =x'$. Conversely, if $x:X\to Y$ is 
cocartesian, and $t:Y\to Y'$ is a vertical isomorphism, then $x.t$ is 
cocartesian. A vertical arrow is cocartesian iff it is invertible.
\end{prop}

\section{The 2-category $Cat/ \B$} \label{2xx}
Let $\B$ be a category. The objects of $Cat /\B$ are the functors with codomain $\B$, like 
the $\pi: \X \to \B$ considered in Section \ref{1xx}. The morphisms are the strictly commutative 
triangles (``functors {\em over} $\B$'')
\begin{equation}\label{triax}\begin{diagram}\X &&\rTo ^{F}&&\Y \\
&\rdTo_{\pi}&&\ldTo_{\pi '}&\\
&&\B,&&
\end{diagram}\end{equation} thus $Cat/\B$ is a standard slice category. 
Note that since the triangle commutes (strictly), $F$ preserves the 
property for arrows of ``being over $a$'', and in particular 
preserves the property of being vertical. To make the slice category 
$Cat/\B$ into a (strict) 2-category, we describe what are the 2-cells:

Given two 
parallel
morphisms in $Cat /\B$, as displayed in
$$
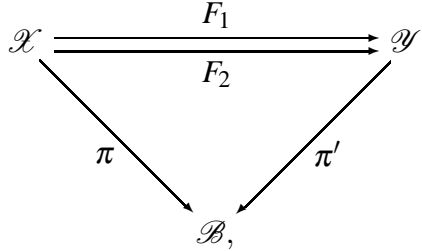
\begin{diagram}\X &&\pile{\rTo^{F_{1}}\\ \rTo_{F_{2}}}&&\Y \\
&\rdTo_{\pi}&&\ldTo_{\pi '}&\\
&&\B,&&
\end{diagram}$$
 the 2-cells between these are taken to be the natural 
transformations $\tau : F_{1} \Rightarrow F_{2}$ which are {\em vertical}, meaning 
that for all objects $X \in \X$, $\tau _{X}$ is a vertical arrow in $\Y$.
Composition (horizontal as well as vertical) is inherited from the 
standard composition of natural transformations in $Cat$. -- Having a 2-category, one has a notion of {\em adjoint} arrows 
``$F\dashv G$ by virtue of 2-cells $\eta $, $\epsilon$''. In particular, 
for the 2-category $Cat/\B$, adjointness of two arrows (functors over 
$\B$) 
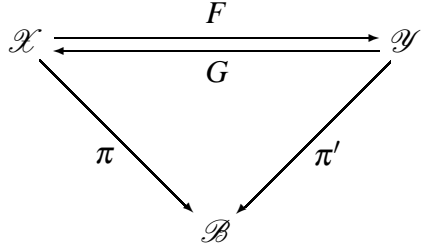
$$\begin{diagram}\X &&\pile{\rTo^{F}\\ \lTo_{G}}&&\Y \\
&\rdTo_{\pi}&&\ldTo_{\pi '}&\\
&&\B&&
\end{diagram}$$
amounts to an ordinary adjointness $\eta ,\epsilon$ between the 
functors $F$ and $G$, subject to the further requirement that $\eta$ 
and $\epsilon$ are {\em vertical} natural transformations. In this 
case, one may write $F\dashv _{\B}G$. If $F\dashv_{\B}G$, 
 then the 
bijection, due to $F\dashv G$, between $\hom  (F(X),Y)$ and $ 
\hom  (X, 
G(Y))$, restricts to a bijection
$$\hom _{a}(F(X),Y)\cong \hom _{a}(X, 
G(Y))$$
whenever $a :A \to B$, $X\in \X_{A}$, $Y\in \Y_{B}$; for, let 
$f:F(X) \to Y$ be an arrow in $\Y$ over $a$. The arrow $X\to G(Y)$ 
corresponding to it under the adjointness is the composite
$\eta_{X}. G(f)$, which an arrow over $a$ since $G(f)$ is so, and 
since $\eta _{X}$ is vertical, by assumption on $\eta$. Similarly, 
verticality of $\epsilon _{Y}$ proves that the inverse correspondence 
preserves the property of being over $a$. 
\begin{prop}[Key Lemma]\label{preservation} Let $F$ and $G$ be vertically adjoint, $F\dashv_{\B}G$, 
as above. Then 
$F$ preserves the property of being cocartesian.\end{prop}
{\bf Proof.} Let $x:X\to Y$ in $\X$ be cocartesian over $a:A\to 
B$, and let $b :B\to C$ be an arbitrary arrow in $\B$. For any $Z\in 
\X _{C}$, the standard naturality 
square for hom set bijections 
induced by the ordinary adjointness  $F\dashv G$ restricts to a commutative square
$$\begin{diagram}\hom_{b}(Y,G(Z)) & \rTo^{\cong} & 
\hom_{b}(F(Y),Z) \\
\dTo ^{x^{*}} && \dTo _{F(x)^{*}}\\
\hom_{a.b}(X, G(Z)) &\rTo_{\cong} 
&\hom_{a.b}(F(X),Z).
\end{diagram}$$
The hom set bijections are displayed horizontally. 
The left hand vertical map is a bijection since $x$ is cocartesian. 
Hence so is the right hand vertical map;  so $F(x)$ is cocartesian.

\section{Opfibrations; cleavages and splittings}\label{3xx}

\begin{defn}Given a functor $\X \to \B$. It is called an {\em 
opfibration} if it has enough cocartesian arrows, in the sense that 
for any arrow $a$ in $\B$ and any $X\in \X$ with $\pi (X) = 
d_{0}(a)$, there 
exists a cocartesian arrow over $a$ with domain $X$.
\end{defn}

A cocartesian arrow $x$ over $a$ is called an ``cocartesian 
{\em lift}'' of $a$. It is a cocartesian  lift of $a$  {\em from $X$} if furthermore 
its domain is $X$.

\medskip

For the remainder of the present Section, $\pi : \X \to \B$ is assumed to 
be an opfibration.

\begin{prop}\label{factx}Given an arrow $z$ in $\X$ with $\pi (z)=a.b$ for 
arrows $a$ and $b$ in $\B$. Then $z$ may be factorized $x.y$ with 
$\pi (x)=a$ and $\pi (y)=b$, with $x$ cocartesian. This facorization 
is unique modulo a unique vertical isomorphism in the middle. And $y$ is 
cocartesian iff $z$ is cocartesian.
\end{prop}
{\bf Proof.} Let $x$ be a cocartesian lift of $a$ with same domain 
as $z$, and construct $y$ over $b$, with $x.y =z$,  using the universal property 
of $x$.  -- The last assertion now follows from Proposition \ref{transx}.

\medskip

A special case is (take $b$ to be the relevant identity arrow):
\begin{prop}\label{fac2} Every arrow $z$ in $\X$ may be factored into 
a cocartesian followed by a vertical arrow. This factorization 
is unique modulo a unique vertical isomorphism.
\end{prop}

\begin{defn} 
A {\em cleavage} for $\X \to \B$ 
consists in a {\em choice}  $X\coc a$ of  a cocartesian lift of 
$a$ from $X$, for 
every $X$ and for every $a $ with $\pi (X)= d_{0}(a )$.
\end{defn}

If there is given a cleavage $\coc$, it is convenient to have a 
separate notation for the codomain of the chosen cocartesian arrow 
$X\coc a$; common notations are $a_{!}(X)$, $\Sigma_{a}(X)$ or $\exists 
_{a}(X)$, or, the one we shall use,
$a_{*}(X):= d_{1}(X \coc a )$,
\begin{equation}\label{cocx}\begin{diagram}X&\rTo ^{X\coc a}& a_{*}(X)
\end{diagram}.\end{equation}

If a cleavage is given, the uniqueness assertions ``modulo vertical 
isomorphisms'' in the Propositions 
\ref{factx} and \ref{fac2} may be sharpened to strict uniqueness, by 
requiring the cocartesian arrows to be provided by the cleavage, in 
an evident way, thus for Proposition \ref{fac2}:  $z: X\to Y$ factors uniquely 
as $X \coc \pi (z)$ followed by a vertical arrow.

There are dual notions: cartesian arrows, cartesian lifts, 
fibrations, with associated cleavage/splitting terminology.
By experience, they are more important than opfibrations. The reason
we discuss opfibrations rather than fibrations is that they are more 
straightforward in so far as variance is concerned. Otherwise, the 
mathematics is the same. Let us right away describe our notation, 
corresponding to (\ref{cocx}), for 
a cleavage of a fibration:
\begin{equation}\label{dcocx}\begin{diagram}
a^{*}(X)& \rTo^{a \dcoc X}& X
\end{diagram}\end{equation}
where now $\pi (X) = d_{1}(a)$. (In Giraud's notation: $\begin{diagram} 
X^{a}&\rTo^{X_{a}}&X\end{diagram}$.)

Consider a functor $F: \X \to \Y$ over $\B$, as displayed in 
(\ref{triax}). Assume both $\X \to \B$ and $\Y \to \B$ are 
opfibrations. Then the $F$ is 
called a {\em morphism of opfibrations} if it takes cocartesian 
arrows to cocartesian arrows. 

Note that both being an opfibration and being a morphism of 
opfibrations are {\em properties} of categories (resp.\ functors) 
over $\B$, not a {\em structure} that is given (resp.\ is 
preserved). In contrast, a cleavage is a structure, and morphisms of 
opfibrations may 
or may not preserve cleavages. 

If the chosen lifts of identity arrows 
are identity arrows, the cleavage is called {\em normalized}; thus, if 
$\pi (X)=A$, we have
$$X\coc 1_{A}=1_{X}.$$ 
If 
furthermore the composite of two chosen cocartesian arrows is again 
chosen, the clevage is called a {\em splitting}:
\begin{equation}\label{splitx}(X\coc a ).(a_{*}(X)\coc b )= X \coc (a.b
).\end{equation} 
Taking the codomain of the two sides of this equation gives in 
particular that
\begin{equation}\label{yXXX}b_{*}(a_{*}(X)) = (a.b)_{*}(X).
\end{equation}

Some of the  literature on (op-) fibrations formulate the theory 
of (op-) fibrations mainly in terms of a cleavages/splittings -- 
which in turn can be reformulated in terms of $Cat$-valued 
pseudofunctors/functors.

\medskip

\noindent {\bf Example.} A group may be considered as a category, in which there is only one 
object, and where all arrows are invertible. If $\X$ and $\B$ are groups, a 
functor $\pi : \X \to \B$ is the same as a group homomorphism. The 
vertical arrows form the kernel of $\pi$. All arrows in $\X$ are 
cocartesian. A group homomorphism 
$\pi :\X \to \B$ is an opfibration iff it is surjective. A 
(normalized) 
cleavage is a (set theoretical) section $s$ of $\pi$ (taking the identity arrow of $\B$ to 
the identity arrow of $\X$). Then $s$ is a splitting iff $s$ is a 
group homomorphism. Not every surjective group homomorphism admits 
such a splitting. So there are opfibrations which do not admit 
splittings.

\medskip

\noindent{\bf Remark.} Any opfibration may, by the axiom of choice, be supplied with a 
normalized cleavage; but not necessarily with a splitting, as the 
example shows.
 However, every 
opfibered category $\X \to \B$ is {\em equivalent} in the 2-category 
$Cat/\B$ to (the underlying opfibration of) a split opfibration over 
$\B$, see Section \ref{sevenx} below. 

This does not contradict the  fact just mentioned about non-existence 
of splittings
of surjective group homomorphism. For, a category equivalent to a 
group need not be a group; it may have several objects.

\section{The ``opfibration'' monad on $Cat/\B$}\label{4xx}
Given an object in $Cat/\B$, i.e.\ a functor $\pi: \X 
\to \B$. One derives from this a new 
$T(\pi: \X \to \B) \in Cat/\B$: it is the comma category $\pi \downarrow 
\B$, equipped with the ``codomain'' functor  to $\B$. Recall that an 
object of $\pi \downarrow 
\B$ is a pair $(X, a )$, where $X\in \X$ and $a $ is an arrow in $\B$ 
with domain $\pi (X)$. The codomain functor $d_{1}:\pi 
\downarrow \B \to \B$ takes this 
object to the codomain of $a$. The arrows in $Cat/\B$ are given by 
arrows in $\X$ together with suitable commutative squares in $\B$; 
for a display, see e.g.\ (\ref{TX-arry}) below.

A more succinct classical description, cf.\ [Gray  1966], is that $T(\pi: \X \to \B)$ is the left 
hand column in the following diagram, in which the square is a 
(strict) pull back:
\begin{equation}\label{pb-dia}\begin{diagram}T (\X \to \B  \SEpbk ) &\rTo & \X \\
\dTo&&\dTo_{\pi}\\
\B ^{{\bf 2}}&\rTo_{d_{0}}&\B\\
\dTo^{d_{1}}&&&\\
\B. &&&
\end{diagram}\end{equation}
Here, $\B^{\bf {2}}$ denotes the standard category of arrows in $\B$, and 
$d_{0}$ and $d_{1}$ are the ``domain'' and ``codomain'' functors, 
respectively.

Thus,  an object $(X, a )$ over $A\in \B$  may be depicted
\begin{equation}\label{TX-obj}\begin{diagram}X&&\\
:&&\\
A&\rTo _{a  }& B
\end{diagram}\end{equation}
(with $X\in \X$, and $\pi (X)=A$).
If $(X',a ')$ is another such object (where $a': A' \to B')$, and if $\beta :B\to B'$ is an 
arrow in $\B$, then
\begin{equation}\label{eqn}\hom _{\beta}^{T(\X)}((X,a  ),(X',a  
'))=
\{ x:X\to X'\mid  
\pi (x).a ' = a .\beta \}\end{equation}
It is a subset of $\hom^{\X}(X,X')$.
Thus, an arrow  ${\bf a}: (X,a  )\to (X',a ')$ over $\beta :B\to 
B'$ may 
be depicted
\begin{equation}\label{TX-arry}\begin{diagram}X&&&&&&\\
&\rdTo^{x}&&&&&&\\
:&&X'&&&&& \\
A&&\rTo^{a}&&B&&&\\
&\rdTo_{\pi (x) }&&&&\rdTo^{\beta }&\\
&&A'&&\rTo_{a'}&&B'
\end{diagram}\end{equation}
(with the bottom square commutative), and it may be denoted $(x,\beta):(X,a) 
\to (X', a')$; it is an arrow over $\beta$.

\medskip

The following shows that $T(\X \to \B)$, as a category over $\B$, has  some canonical 
cocartesian arrows: 

\begin{prop} Given an object $(X, a )$ in $T(\X \to \B)$  over $B$, and 
given an arrow $b:B\to C$ in $\B$. Then there is a canonical cocartesian arrow over $b$ from $(X, 
a)$ to $(X, a .b)$, depicted in

\begin{equation}\label{zuy}\begin{diagram}X&&&&&&&\\
&\rdTo^{1_{X}}&&&&& &\\
:&&X&&&&&& \\
A&&\rTo^{a}&&B&&&\\
&\rdTo_{1_{A} }&&&&\rdTo^{b }&\\
&&A&&\rTo_{a.b}&&C.
\end{diagram}\end{equation}
\end{prop}
We denote this arrow in $T(\X \to \B )$ by $((X,a);b)$, it is an arrow over $b$  from 
$(X,a)$  to $(X,a.b)$. We may define a cleavage by putting $(X,a)\coc 
b):= ((X,a);b$, thus $$\begin{diagram} (X,a) & 
\rTo^{(X,a)\coc b}&  
(X,a.b).\end{diagram}$$

\medskip

\begin{sloppypar}\noindent{\bf Proof.} To see that the depicted arrow $((X,a);b):(X,a) 
\to (X,a.b)$ is cocartesian over $b:B \to C$, let 
$c:C\to D$ and let $(Z, \delta )$ be an object over $D$ (where 
$\delta :D' \to D$ and $\pi(Z) = D'$). Then as subsets of 
$\hom(X,Z)$, we see that
$\hom _{c}((X, a.b), (Z, \delta ))$ consists of those $h:X\to Z$ which satisfy 
$\pi (h). \delta = (a.b).c$, and $\hom _{b.c}((X,a),(Z,\delta))$ of those 
$h:X\to Z$ which satisfy $\pi (h). \delta = a .(b.c)$, and these 
two subsets are equal. Precomposition with (\ref{zuy}) is provided by 
precomposition with the $\X$ component which here is $1_{X}$, so is 
indeed the identity mapping of the described subset onto 
itself.\end{sloppypar}

\medskip  

From the Proposition immediately follows that $T(\X \to \B) \to \B$ has 
sufficiently many cocartesian arrows to deserve the title 
``opfibration'', and in fact, the very construction of the canonical 
cocartesian arrows shows that it provides this opfibration with a 
{\em splitting}: the codomain of $((X,a);b)$ is $(X,a.b)$, and clearly
$$((X,a);b).((X,a.b);c)= ((X,a);b.c).$$

\medskip 

The functorial character of $T$ is straightforward from the 
construction (\ref{pb-dia}). Explicitly: for a 
functor over $\B$, as depicted in (\ref{triax}), $T(F)$ is the 
functor $T(\X \to \B) \to T(\Y \to \B)$ which on objects is 
given by $T(F)(X,\alpha ) = (F(X),\alpha )$,
and on arrows $T(F)(x)= F(x)$ (if $x$ satisfies the equation 
(\ref{eqn}), then so does $F(x)$ (with $\pi$ replaced by $\pi'$)).

We make $T$ into a monad by supplying natural transformations $y:Id 
\Rightarrow T$ and $m: T^{2}\Rightarrow T$. 
We describe first the objects of 
$T^{n}(\X \to \B)$: an object  over $A\in \B$ of $T^{n}(\X \to \B$) 
may  be depicted
$$\begin{diagram}
X&&&&&&\\
:&&&&&&\\
A^{(n)}&\rTo&\ldots&\rTo &A''&\rTo &A'&\rTo&A .
\end{diagram}$$
The unit $y$ and multiplication $m$ of the monad  come 
about from the units and the composition in the chain of $A^{(i)}$s 
(and this makes  the unit and associative laws for $y$ and $m$ evident).
Thus (if we suppress $\pi$ from notation), the functor $y_{\X}:\X \to T(\X )$ 
takes an $X\in \X_{A}$ to the configuration
\begin{equation*}\begin{diagram}X&&\\
:&&\\
A&\rTo_{1_{A}} & A,
\end{diagram}
\end{equation*}
or $X\mapsto (X, 1_{A})$. 
Similarly, the object $(X,a,b)\in T^{2}(\X)$ 
goes by $\mu _{\X}$ to the object $(X, a.b) 
\in T(\X)$. 

Consider an object $(X,a )\in T(\X)$, as depicted in (\ref{TX-obj});
then $T(y_{\X })(X, a )$ is the object $(X,1,a )$, depicted 
in

\begin{equation}\label{Tyy}\begin{diagram}
X&&&&&\\
:&&&&&\\
A&\rTo_{1_{A}}&A&\rTo_{a}&B,
\end{diagram}
\end{equation}

and $y_{T(\X)}(X, a)$ is the object $(X, a , 1)$ depicted in

\begin{equation}\label{yTy}\begin{diagram}
X&&&&&\\
:&&&&&\\
A&\rTo_{a}&B&\rTo_{1_{B}}&B,
\end{diagram}
\end{equation}

The reader may, as an exercise, describe a vertical arrow in $T^{2}(\X )$ from 
the first of these objects to the second.

The endofunctor $T$ on $Cat/\B$ is clearly  canonically enriched over $Cat$: 
its value on a 2-cell between 1-cells $F$ and $G$, i.e.\ on a vertical 
natural transformation $t:F\Rightarrow G$, is the vertical natural 
transformation whose instantiation at $(X,a ) \in T(\X)$ (where 
$a:A\to B$) 
associates   $t_{X}:F(X) \to G(X)$ (it may be viewed as an element in 
$\hom^{\Y}_{A}((F(X),a )(G(X),a ))$ since it satisfies an equation 
like (\ref{eqn}): $\pi '(t_{X}).a = a .1_{A}$, because $\pi 
'(t_{X})$ is an identity arrow). And $y$ and $m$ are 2-natural.

\medskip

The  name ``opfibration monad'' should really be, more precisely: the 
``split-opfibration monad''; we have already argued that every 
$T(\X \to \B )$ has canonically  the structure of  a {\em split} opfibration; 
we shall see that $T(\X \to \B)$ is 
the free such on $\X \to \B$, in fact, we shall see that the category of 
Eilenberg-Moore algebras for this monad is the category of split 
opfibrations over $\B$.

\section{KZ aspects}\label{KZA}
The monad $(T,y,m)$ on $Cat/\B$ has now been enhanced to a (strict) 
2-monad. It can be even further enhanced, namely to a KZ monad in the 
sense of [Kock 1973, 1995]. This 
means  that it is provided with a modification
$\lambda : T(y) \Rightarrow y_{T}$, with certain equational 
properties. 
Concretely, it here means that for each object $\X= (\pi: \X\to \B)$ 
in $Cat/\B$, there is provided a vertical natural transformation $\lambda = 
\lambda_{\X}$,
$$\begin{diagram}T(\X)&\pile{\rTo^{T(y_{\X})} \\ \Downarrow\lambda \\ 
\rTo_{y_{T(\X)}}}& T^{2} (\X) 
\end{diagram}$$
satisfying two ``whiskering'' equations (\ref{whisy}) and (\ref{whis1}) below.
Consider an object $(X,a )\in T(\X)$, as in (\ref{TX-obj}); we describe 
$\lambda _{(X,a)}$. Recall the description of the objects $T(y_{\X})(X,a)$ and 
$y_{T(\X)}(X,a )$ given in (\ref{Tyy}) and (\ref{yTy}), 
respectively.
Then $\lambda _{(X,a)}$ is the arrow in $T^{2}(\X)$ given by the 
following. The unnamed arrows are identity arrows.

\begin{equation}\label{lambday}\begin{diagram}X&&&&&&&&&&\\
&\rdTo^{}&&&&& &&&&\\
:&&X&&&&&& &&&\\
A&\rTo^{}&&A&\rTo^{a} &&B\\
&\rdTo_{ }&&&\rdTo^{a}&&&\rdTo^{}&\\
&&A&\rTo_{a}&&B&\rTo_{}&&B&
\end{diagram}\end{equation}
(Cf.\ the Exercise comparing (\ref{Tyy}) and (\ref{yTy})).
It is clear that if $a$ is an identity arrow in $\B$, then $\lambda 
_{(X,a)}$ is an  identity arrow in $T(\X)$, which  
implies that the whiskering 
\begin{equation}\label{whisy}\begin{diagram}\X&\rTo^{y_{\X}}&T(\X)&\pile{\rTo^{T(y_{X})} \\ \Downarrow\lambda \\ 
\rTo_{y_{T(\X)}}}&T^{2}(\X)
\end{diagram}\end{equation}
is an identity natural transformation. Also, applying $m_{\X}$ to 
the arrow in $T^{2}(\X)$ exhibited in (\ref{lambday}) produces the identity 
arrow of $(X, a )$ in $T(\X )$, and this implies that the whiskering
\begin{equation}\label{whis1}\begin{diagram}T(\X)&\pile{\rTo^{T(y_{X}) }\\ \Downarrow\lambda \\ 
\rTo_{y_{T(\X)}}}&T^{2}(\X)&\rTo^{m_{\X}}&T(\X)
\end{diagram}\end{equation}
is the identity natural transformation. These two whiskering 
equations are what a modification $\lambda$ should satisfy in order to 
make a (strict) monad on a 2-category into a KZ monad, cf.\ [Kock 1995], Axioms T0-T3 
(with T0 and T3 being redundant if $(T,y,m)$ is a strict monad, which 
is the case here).

Whenever one has a monad $T$ on a category ${\mathfrak C}$, one has the category 
of (Eilenberg -Moore) algebras for 
it, i.e.\ an object $\X \in {\mathfrak C}$ together with a ``structure'' map $\xi 
:T(\X ) \to \X $, satisfying the standard unit- and associativity equations.

Recall  that in a 2-category, the notion of adjointness between 1-cells 
makes sense.

 We denote objects in a 2-category ${\mathfrak C}$ by script 
letters like $\X$, because of the example we have in mind. Also, we 
compose the arrows in ${\mathfrak C}$ from right to left.

Here is a basic construction in the context of KZ monads (cf.\ [Kock 
1973, 1995]).
To produce an adjointness $\xi \dashv 
y_{\X}$ out of an Eilenberg-Moore algebra $\X , \xi$, 
we produce unit $\eta$ and counit $\epsilon$. For the counit $\epsilon$, 
we just take the identity 2-cell on $\xi \circ y_{\X}. = id_{\X}$. The unit 
$\eta$ is constructed using $\lambda$: we have the 2-cell obtained by 
whiskering $\lambda _{\X}$ with $T(\xi )$:
\begin{equation}\label{whisk}\begin{diagram}T(\X)&\pile{\rTo^{T(y_{\X}) }\\ \Downarrow\lambda \\ 
\rTo_{y_{T(\X)}}}&T^{2}(\X)&\rTo^{T(\xi )}&T(\X)
\end{diagram}.\end{equation}
The top composite is an identity 1-cell, since $\xi \circ y_{\X} = 1_{\X}$. 
The lower composite may be rewritten as $ y_{\X}\circ \xi$, using 
naturality of $y$ w.r.to $\xi$. So the whiskering (\ref{whis}) gives 
a 2-cell
\begin{equation}\label{whis}\begin{diagram}T(\X) & \pile{\rTo^{id}\\ 
\Downarrow \\ \rTo_{y_{\X}\circ 
\xi}}&T(\X)
\end{diagram},\end{equation}
 and this is the unit $\eta$ of the adjointness 
$\xi \dashv y_{\X}$. The triangle equation holds by virtue one of the 
whiskering equation (\ref{whisy}). In particular, since $m_{\X}$ is a 
$T$-homomorphism, we have $m_{\X}\dashv y_{T(\X )}$.

\medskip

If $\xi \dashv y_{\X}$,
 it does 
not conversely follow that $\xi$ is an Eilenberg-Moore algebra 
structure,
 since the associative law $\xi \circ T(\xi ) = \xi \circ m$ may not 
hold strictly; but we have from [Kock 1973, 1995]
\begin{prop}If $T$ is a KZ monad, and if $\xi: T(\X) \to \X$ is a left 
adjoint for $y_{\X}$, then there is a canonical isomorphism (2-cell) 
$\alpha$ between 
the two 1-cells  $\xi \circ T(\xi)$ and $\xi \circ m_{\X}$.
\end{prop}
{\bf Proof.} Since $\xi \dashv y_{\X}$, and $T$ is a 
2-functor, it follows that $T(\xi ) \dashv T(y_{\X})$. Also, since 
$m_{\X}:T^{2}(\X ) \to T(\X )$ is an Eilenberg-Moore structure by 
general monad theory, and because $T$ is KZ, it follows that 
$m_{\X}\dashv y_{T(X)}$. We therefore have that
$$\xi \circ T(\xi ) \dashv T(y_{\X}) \circ y_{\X}$$
and partly that
$$\xi \circ m_{\X} \dashv y_{T(\X)}\circ y_{\X};$$
but the two right hand sides here are equal, by naturality of $y$, so 
it follows that the two left adjoints exhibited are  
canonically isomorphic.
\medskip

\begin{sloppypar}For a KZ monad, it can be proved (cf.\ loc.cit)  
that the canonical isomorphism $\alpha$ 
satisfies those coherence equations which make $\X, \xi$ into an 
Eilenberg-Moore pseudo-algebra, in the standard sense of 2-dimensional 
monad theory; vice versa, a pseudo-algebra $\X , \xi , \alpha $ in the standard sense
has $\xi \dashv y_{X}$.\end{sloppypar}

\medskip

\section{Cleavages and splittings in terms of the opfibration 
monad}\label{6xx}
We now return to the case where $\C = Cat/\B$, and where $T$ is the 
``opfibration KZ monad'' described. Recall that $\X$ is also used as a 
shorthand for an object $\pi : \X \to \B$ in $Cat/\X$. 
\begin{thm}1) Assume $\xi : T(\X ) \to \X$ is vertically left adjoint to $y_{\X}$ 
in $Cat/\X$,   and with the identity functor 
on $\X$ as counit (so $\xi\circ y_{\X} =id_{\X}$).\footnote{[Gray 
1966] 
introduced  the short hand ``lali'' for a left adjoint left 
inverse, like $\xi$.} Then $\X$ carries 
a canonical structure of  normalized cleavage. Conversely, a 
normalized cleavage defines a functor $\xi $, left adjoint to $y_{\X}$ 
and with $\xi\circ y_{\X} = id$.

2) This normalized cleavage is a splitting if and only if $\xi$ is strictly 
associative, $T(\xi )\circ \xi = m\circ \xi$.

3) If $(\X , \xi )$ and $(\X ',\xi ')$ are strict $T$-algebras, then a 
functor $\X \to \X'$ over $\B$ is a $T$-homomorphism iff it preserves 
the corresponding splittings strictly.
\end{thm}
{\bf Proof/Construction.} Given $\xi$. Let $X\in \X_{A}$, and let $a:A \to B$ in 
$\B$. We produce our candidate for a cocartesian lift $X\coc a: X\to 
a_{*}(X)$  by applying $\xi$ to 
the arrow $((X,1_{A}); a)$ in $T(\X)$ exhibited in the following 

\begin{equation}\label{zxxy}\begin{diagram}X&&&&&&&\\
&\rdTo^{1_{X}}&&&&& &\\
:&&X&&&&& \\
A&&\rTo^{1_{A}}&A&&\\
&\rdTo_{1_{A} }&&&\rdTo^{ a}&\\
&&A&&\rTo_{a}&B
\end{diagram}\end{equation}
Thus $a_{*}(X)$ is $\xi (X,a)$. Note that (\ref{zxxy}) is a special case of the general canonical 
cocartesian arrow (\ref{zuy})  in $T(\X)$. So (\ref{zxxy}) is 
cocartesian, and therefore, applying $\xi$ to it gives, by the Key 
Lemma  
(Proposition \ref{preservation}), a cocartesian 
arrow in $\X$. Its domain is $X$, since the top line in the 
(\ref{zxxy}) 
is $y_{\X}(X)$ and $\xi \circ y_{\X} = id_{\X}$. It lives over $a$, 
since the right hand slanted arrow in ({\ref{zxxy}) is $a$. 
Thus we have constructed a cleavage for $\X$. The fact 
that the constructed cleavage is normalized follows because $\xi 
\circ y_{\X}$ is the identity functor, and  because $\xi$, as 
a functor, takes identity arrows  to  identity arrows. 

Conversely, given a normalized cleavage $\coc$ of $\X \to \B$. 
Then a functor $\xi : T(\X)\to 
\X$ is constructed as follows:  on objects, we put
$\xi (X,a ):=a_{*}(X)$ (= $ d_{1}(X\coc a ))$, and on 
morphism we use the universal property of cocartesian arrows. More 
explicitly, $\xi$ applied to the arrow ${\bf a}$ in $T(\X)$ over 
$\beta$ displayed in 
(\ref{TX-arry}) is the unique arrow $\xi({\bf a})$ over $\beta $ which makes the square in $\X$
\begin{equation}\label{xi-arr}\begin{diagram}
X&\rTo^{X\coc a}&a_{*}(X)\\
\dTo^{x}&&\dDashto _{\xi ({\bf a})}\\
X'&\rTo_{X'\coc a '}&a_{*}'(X')\\
\end{diagram}\end{equation} 
commute.

The  fact that $\xi \circ y_{\X} 
= id$ to left to the reader. The fact that $\xi$ is indeed a 
left adjoint for $y_{\X}$ follows because it solves a universal 
problem; the unit of the adjunction at the object $(X,a)$ is the arrow  $(X,a) \to 
(a_{*}X, 1_{B})$ in $T(\X )$  given by the diagram

\begin{equation*}\begin{diagram}
X&&&&&&&\\
&\rdTo^{X\coc a}&&&&& &\\
:&&a_{*}(X)&&&&& \\
A&&\rTo^{a}&B&&\\
&\rdTo_{a }&&&\rdTo^{ 1_{B}}&\\
&&B&&\rTo_{1_{B}}&B.
\end{diagram}
\end{equation*}
This proves the assertion 1) of the Theorem.

Assume next that $\xi$ is associative, i.e.\ $\xi \circ T(\xi ) = \xi 
\circ m$. The 
arrows picked out by the clevage $\coc$ derived from $\xi$   are those that are of the 
form: $\xi$ applied to an arrow in $T(\X )$ of the form $((X,1_{A}); a)$, as in 
(\ref{zxxy}); whereas $\xi$ applied to a more general canonical cocartesian 
arrow in $T(\X )$ of the form $((X,a);b)$, as exhibited in (\ref{zuy}),
is 
not apriori picked out by the cleavage $\coc$. However, we have, with 
$\coc$ and the resulting $a_{*}(X)$ derived from $\xi$, the following

\begin{lemma}Assume that $\xi$  is  associative. Then
$$\xi ((X,a);b) = a_{*}(X) \coc b.$$
\end{lemma}
{\bf Proof.} Consider the following arrow over $b$ in $T^{2}(\X)$:

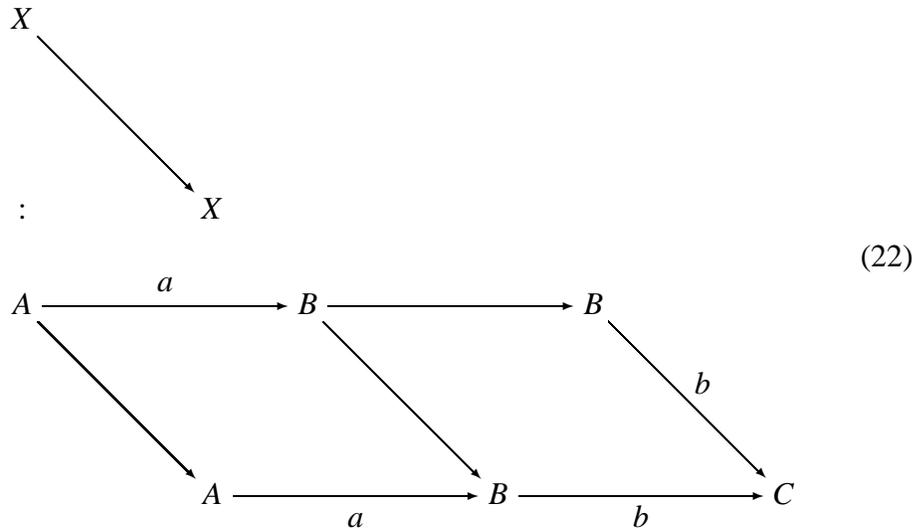
\begin{equation}\label{onexy}\begin{diagram}X&&&&&&&&&&\\
&\rdTo^{}&&&&& &&&&\\
:&&X&&&&&& &&&\\
A&\rTo^{a}&&B&\rTo &&B\\
&\rdTo_{ }&&&\rdTo^{}&&&\rdTo^{b}&\\
&&A&\rTo_{a}&&B&\rTo_{b}&&C&
\end{diagram}\end{equation}
(unnamed arrows are identity arrows); applying $m_{\X}$ yields

\begin{equation}\label{twoxy}\begin{diagram}
X&&&&&&&\\
&\rdTo&&&&& &\\
:&&X&&&&& \\
A&&\rTo^{a}&B&&\\
&\rdTo&&&\rdTo^{ b}&\\
&&A&&\rTo_{a.b}&C
\end{diagram}
\end{equation}

i.e.\ $((X,a);b)$; whereas applying $T(\xi )$ 
to (\ref{onexy}) 
yields

\begin{equation}\label{threexy}\begin{diagram}
a_{*}(X)&&&&&&&\\
&\rdTo^{}&&&&& &\\
:&&a_{*}(X)&&&&& \\
B&&\rTo^{}&B&&\\
&\rdTo&&&\rdTo^{ b}&\\
&&B&&\rTo_{b}&C
\end{diagram}
\end{equation}
(since $\xi (X,a)=a_{*}(X)$ and since the first square of (\ref{onexy}) is an 
identity arrow in $T(\X)$). By the strict associativity of $\xi$, the 
value of $\xi$ on (\ref{twoxy}) and (\ref{threexy}) is the same, and 
these values  are $\xi ((X,a);b)$ and  $a_{*}(X)\coc b$, respectively. This proves the Lemma.

 To prove the splitting 
condition (\ref{splitx}), let an arrow 
$b:B\to C$ be given. 
Consider the composite in $T(\X)$
\begin{equation}\label{compxx}\begin{diagram}(X,1_{A})&\rTo^{((X,1_{A});a)}&(X,a)&\rTo^{((X,a);b)}&(X,a.b).
\end{diagram}\end{equation}
 One sees that applying $\xi$ (using the Lemma for the second factor)
 gives the composite
composite $(X\coc a). (a_{*}(X)\coc b)$.  On the other hand, the 
composite (\ref{compxx}) is $((X,1_{A}); a.b)$, which by $\xi$ gives 
$X\coc (a.b)$. This proves that the cleavage $\coc$ produced by a strictly associative $\xi$ 
 is in fact a splitting.

Conversely, given a splitting $\coc$, then since $\coc$ is in 
particular a normalized cleavage, it gives rise to a functor 
$\xi :T(\X) \to \X$, with $\xi \circ y_{\X}$ the identity functor on 
$\X$, as described above.
It remains to prove that $\xi$ satisfies   the associative law $\xi \circ T(\xi 
)= \xi \circ m :T^{2}(\X ) \to \X$. 
Consider an object $(X,a , b)$ in $T^{2}(\X 
)_{B}$, as displayed in
$$\begin{diagram}X&&&&\\
:&&&&\\
A&\rTo_{a}&B&\rTo_{b}&C.
\end{diagram}$$
Then
$$\xi (T(\xi)(X,a ,b)) = \xi (a_{*}(X), b) 
= b _{*}(a_{*}(X)).$$ On the other hand
$$\xi (m (X,a , b))= \xi (X, a. b) = 
(a.b)_{*}(X),$$ and then 
(\ref{yXXX})  gives the 
associativity result, in so far as objects of $T^{2}(\X)$ is 
concerned. Next, consider a morphism in $T^{2}(\X )$ from
$(X,a ,b )$ to $(X',a',b')$ 
 given by $(x,\beta,\gamma
)$ (with $\alpha= \pi (x)$), displayed as the full arrows in the 
three-dimensional diagram 
$$\begin{diagram}
X&&\rDashto^{X\coc a} &&a_{*}(X)&&\rDashto^{a_{*}(X)\coc b}&&b_{*}(a_{*}(X)) &&\\
&\rdTo^{x}&&&&\rdDashto^{\xi (x,\beta)}&&&&\rdDashto^{\xi (\xi 
(x,\beta ),\gamma )}&\\
:&&X'&&\rDashto^{X'\coc a'}&&a'_{*}(X') &&\rDashto^{a_{*}'(X')\coc 
b'}&&b'_{*}(a'_{*}(X'))\\
A&&\rTo^{a} &&B&&\rTo^{b}&&C&&:\\
&\rdTo^{\alpha}&&&&\rdTo^{\beta }&&&&\rdTo^{\gamma}&\\
&&A'&&\rTo_{a'} &&B'&&\rTo_{b'}&&C'\\
\end{diagram}$$
 The  
slanted dotted
arrows in the top layer are, 
by construction of the value of 
$\xi$ on arrows in $T(\X )$, the unique ones (over $\beta$ and $\gamma$, 
respectively) which make the squares on the top commute. 

So $\xi 
\circ T(\xi )$ applied to the given arrow in $T^{2}(\X )$ is the 
rightmost slanted arrow on the top. On the other hand, $ m_{\X}$ 
applied to the given morphism in $T^{2}(\X )$ is given by the 
full arrows in 
\begin{equation}\label{three}\begin{diagram}X&&\rDashto^{X\coc 
(a.b)}&&(a.b)_{*}(X)&\\
&\rdTo^{x}&&&&\rdDashto ^{\xi (\xi(x, \beta ),\gamma)} &\\
:&&X'&\rDashto^{X'\coc (a'.b')}&&&(a'.b')_{*}(X') \\
A&&\rTo^{a.b}&&C&&:\\
&\rdTo_{\alpha }&&&&\rdTo^{\gamma }&\\
&&A'&&\rTo_{a'.b'}&&C'
\end{diagram}\end{equation}
so $\xi$ applied to it is  
unique one over $\gamma$ making the top square commute. But now by 
the splitting condition for $\coc$,
equation (\ref{splitx}), we conclude that the composite of $\coc$ arrows in the 
previous diagram equals the $\coc$ arrow in the present one, so by 
uniqueness of chosen cocartesian lifts of $\gamma$, we conclude that the two desired 
arrows in $\X$ agree, proving $\xi \circ m=\xi \circ T(\xi )$.

It is clear that the two processes $\xi \leftrightarrow \coc$ are 
mutually inverse, and so the assertion 2) is proved.

Finally, consider two strict $T$-algebras $(\X ,\xi ) $ and $(\X', \xi 
')$. We need to prove that a functor $F: \X \to \X'$ over $\B$ 
is compatible with $T$-algebra structures $\xi$ and $\xi'$  iff it 
is compatible with the associated splittings $\coc$ and $\coc '$. If 
$F$ is compatible with the algebra structures, we get that $F$  
takes the  cocartesian arrows $X\coc a$ in $\X$ to the 
cocartesian arrow $F(X) \coc ' a$ in $\X'$; this follows by 
considering the canonical cocartesian arrow (\ref{zxxy}) and applying 
the two functors (assumed equal) $F\circ \xi $ and $\xi ' \circ T(F)$ 
to it. Conversely, if $F(X\coc a )= F(X) \coc 'a$,
 the functors $F\circ \xi $ and $\xi ' \circ T(F)$ give equal value 
 in $\X'$ on the object $(X,a) \in T(\X )$; for,
$$F(\xi (X, a ))= F(d_{1}(X\coc a )) = d_{1}(F(X\coc a
)) = d_{1}(F(X) \coc 'a )$$ $$ = \xi' (F(X), a) = \xi' ( 
T(F)(X,a )).$$
The two functors in question then agree also on morphisms; this 
follows from the fact that they agree on objects, and from the fact that their 
values on morphisms is determined by universal properties. This 
proves assertion 3) and thus the Theorem.

\medskip

The two last assertions of the Theorem immediately lead to
\begin{cor} The category of split (op-)fibrations over $\B$, and strict 
splitting preserving functors, is monadic over $Cat /\B$ by the 
KZ monad $T$.
\end{cor}

Also, the category of opfibrations with a cleavage is the category of 
pseudo-algebras for the monad $T$, with morphisms: functors over $\B$  preserving the 
clevages strictly; and with pseudo-morphisms: the functors over $\B$ 
that preserve cocartesian arrows. We rephrain from making these 
``pseudo-'' notions explicit, but it should be mentioned that those natural isomorphisms 
that occur in  the precise definition of the ``pseudo-'' notions 
 for the case of KZ monads automatically satisfy the 
coherence conditions that usually must be required for such 
isomorphisms, because they solve universal problems.

\section{Replacing cleavages with splittings}\label{sevenx}
Every opfibration   admits a normalized clevage (granted the axiom of 
choice); 
but as remarked in the Example and Remark at the end of Section 
\ref{3xx}, it may not admit a splitting. On the other hand, one has 
([Giraud] I.2): 
\begin{thm}Every opfibration $\pi :\X \to \B$ is equivalent (in the 
2-category $Cat/\B$) to one with a splitting.
\end{thm}
{\bf Proof.} Choose a normalized cleavage $\coc$ for $\X \to \B$, and 
construct the corresponding left adjoint left inverse $\xi : T(\X ) \to \X$ for 
$y_{\X}$. Then take the full image of $\xi$  in $\X$. Recall that the full 
image $\F (\xi )$ of a functor $\xi :\Y \to \X$ has the same objects as $\Y$, and 
that the set of arrows  $Y_{1}\to Y_{2}$ in $\F (\xi )$ is the set of 
arrows $\xi (Y_{1}) \to \xi (Y_{2})$ in $\X$. (If you want disjoint 
hom-sets in $\F (\xi )$, then put on some labels $Y_{i}$). 
There is an evident factorization of $\xi$:
$$\begin{diagram}\Y &\rTo ^{\xi_{1}}& \F (\xi ) & \rTo^{\xi_{2}}& \X 
\end{diagram}$$
where $\xi _{1}$ is bijective on objects and $\xi _{2}$ is full and 
faithful.
If $\xi$ is 
surjective on objects, $\xi_{2}$ is an equivalence of categories. 
In our case, $\xi: T(\X) \to \X$ is surjective on objects, because 
 $\xi \circ y_{\X}$ is the identity functor on $\X$. Also, it is easy 
to see that the full image construction $\F (\xi )$ respects the 
``augmentations'' to $\B$. So the proof is completed if we can 
provide $\F (\xi )\to \B$ with the structure of a split opfibration. 
But $T(\X )$ carries a 
canonical splitting, which we can transfer to $\F (\xi )$ using $\xi$.
If we denote the canonical splitting of $T(\X)$ by $\coc_{0}$, then 
we define a cleavage $\coc_{1}$ on $\F (\xi )$ by the formula
$$\begin{diagram}Y\coc_{1}a\quad := \quad  
\xi_{1}(Y&\rTo^{Y\coc_{0}a}&a_{*}(Y))
\end{diagram}$$
where $a_{*}(Y)$ denotes the codomain of $Y\coc_{0}a$.
We note that this arrow in $\F (\xi )$ has indeed the correct domain, 
namely $\xi _{1} (Y) =Y$; we also note that its codomain is again 
$a_{*}(Y)$. To prove that the cleavage $\coc_{1}$ is a 
splitting, we consider  the equation
$(Y\coc_{0}a).(a_{*}Y\coc _{0}b) = Y\coc_{0} (a.b)$, which holds, 
since $\coc_{0}$ is a splitting. We then get the desired equation for 
$\coc_{1}$ 
using that $\xi_{1}$ is a functor, and applying the definition of 
$\coc_{1}$ to each of the three terms of the equation.

\section{Comparisons}\label{9x}
Apparently, Chevalley was the one to formulate the 
notion of fibration in adjointness terms; [Gray 1966] (p.\ 56) uses the term 
``Chevalley Criterion'' for the following (which I here state
 for opfibrations rather than fibrations):

\medskip 

\noindent {\em a functor
  $\pi : \X  \to \B$ is an opfibration iff the canonical functor
$\overline{\pi} :\X^{{\bf 2}}\to \pi \downarrow \B$ admits a left adjoint right inverse 
(lari) $K$.}

\medskip

\noindent Recall that the criterion considered presently is that
$y_{\X}:\X \to T(\X ) $  admits a (vertical) left adjoint {\em left} inverse 
(lali) $\xi$; and recall also that $T(\X )$ is $\pi \downarrow 
\B$ (seen as a category over $\B$). For a lari, it is the counit of the adjunction  which carries the information 
(the unit being an identity); for a lali, it is the unit which 
carries the information. [Street 1974] analyzed (p.118-119), in abstract 
2-categorical terms, that the data of $K$ and $\xi$ are equivalent. I shall 
here describe, in elementary terms, the passage from $\xi$ to $K$, 
and describe the counit for the lari adjunction $K \dashv 
\overline{\pi}$. Recall the notation applied in the present article:
for $(X,a)\in  \pi \downarrow \B$, $\xi$ returns the value 
$a_{*}(X)$, and the unit of the adjunction is essentially $X\coc a$.
 Out of this data, one constructs $K:  \pi \downarrow \B \to \X ^{{\bf 
2}}$ by sending $(X,a)$ to the arrow $X\coc a : X \to a_{*}(X)$. The 
counit of the lari adjunction $K \dashv \overline{\pi}$, instantiated 
at an object $x:X\to Y$ in $\X ^{{\bf 2}}$ is an arrow 
in $\X ^{{\bf 2}}$, namely the commutative square
$$\begin{diagram}X&\rTo^{X\coc \pi (x)}&\pi(x)_{*}(X)\\
\dTo^{1_{X}}&&\dDashto \\
X&\rTo_{x}&Y
\end{diagram}$$
where the right hand vertical arrow comes from the universal property 
of the cocartesian arrow on the top. 

\medskip

Street  was probably the first to observe that opfibrations could 
be described as pseudo-algebras for a KZ  monad;  in fact, in  
[Street 1974] p.\ 118, he uses this 
description as his {\em definition} of the notion of opfibration,  so 
therefore, no proof is given. Also, loc.cit.\ gives no proof  of 
the fact that split opfibrations then are the are the strict 
algebras. So in this sense, Section \ref{6xx} of the present article only supplements 
loc.cit.\ by providing elementary proofs of these facts.

\small
\noindent Diagrams were made with Paul Taylor's package.

\bigskip

\noindent Anders Kock, University of Aarhus, 

\noindent mail: kock (at) math.au.dk

\noindent December 2013.
\end{document}